\documentclass[11pt,amsfonts]{article}
\usepackage{graphicx}
\usepackage{latexsym}
\usepackage{amssymb}
\usepackage{amsmath}
\usepackage{layout}
\newtheorem{prop}{Proposition}
\newtheorem{lemma}{Lemma}

\newtheorem{theorem}{Theorem}
\newtheorem{remark}{Remark}

\def\real{{\mathord{{\rm I\kern-2.8pt R}}}}        % Fake blackboard bold R.
\def\inte{{\mathord{{\rm I\kern-2.8pt N}}}}

\def\sZZ{{\rm Z\kern-2.8ptem{}Z}}

\def\z{{\mathchoice
  {\sZZ}
  {\sZZ}
  {\rm Z\kern-0.30em{}Z}
  {\rm Z\kern-0.25em{}Z} }}
\def\sQQ{{\kern 0.27em \vrule height1.45ex width0.03em depth0em
          \kern-0.30em \rm Q}}
\def\qu{{\mathchoice
    {\sQQ}
    {\sQQ}
  {\kern 0.225em \vrule height1.05ex width0.025em depth0em \kern-0.25em \rm Q}
  {\kern 0.180em \vrule height0.78ex width0.020em depth0em \kern-0.20em \rm Q}
        }}
\def\sCC{{\kern 0.27em \vrule height1.45ex width0.03em depth0em
          \kern-0.30em \rm C}}
\def\complex{{\mathchoice
    {\sCC}
    {\sCC}
  {\kern 0.225em \vrule height1.05ex width0.025em depth0em \kern-0.25em \rm C}
  {\kern 0.180em \vrule height0.78ex width0.020em depth0em \kern-0.20em \rm C}
        }}

\newcommand{\ba}{\begin{array}}
\newcommand{\ea}{\end{array}}
\newcommand{\be}{\begin{equation}}
\newcommand{\ee}{\end{equation}}
\newcommand{\bea}{\begin{eqnarray}}
\newcommand{\eea}{\end{eqnarray}}
\newcommand{\beaa}{\begin{eqnarray*}}
\newcommand{\eeaa}{\end{eqnarray*}}

\newcommand{\eps}{\varepsilon}

\def\b{\beta}

\def\z{\zeta}

\font\tenmath=msbm10 \font\sevenmath=msbm7 \font\fivemath=msbm5
\newfam\mathfam \textfont\mathfam=\tenmath
\scriptfont\mathfam=\sevenmath \scriptscriptfont\mathfam=\fivemath

\def \b{\noindent}

\def \={{\buildrel {\rm (law)} \over =}}

\def\qed{ \hfill \vrule width.25cm height.25cm depth0cm\smallskip}

\newcommand{\basa}{\begin{assumption}}
\newcommand{\easa}{\end{assumption}}

\newcommand{\bas}{\begin{assum}}
\newcommand{\eas}{\end{assum}}

\newcommand{\ignore}[1]{}
\begin{document}

\renewcommand{\thefootnote}{\fnsymbol{footnote}}

\date{ }
\title{Hsu-Robbins and Spitzer's theorems for the variations of fractional Brownian motion}
\author{Ciprian A. Tudor\\
  SAMOS/MATISSE,
Centre d'Economie de La Sorbonne,\\ Universit\'e de
Panth\'eon-Sorbonne Paris 1,\\90, rue de Tolbiac, 75634 Paris Cedex
13, France.\\tudor@univ-paris1.fr\vspace*{0.1in}}
\maketitle

\begin{abstract}
Using recent results on the behavior of multiple Wiener-It\^o integrals based on Stein's method, we prove Hsu-Robbins and Spitzer's theorems for sequences of correlated random variables related to the increments of the fractional Brownian motion.

\end{abstract}

\vskip0.5cm

{\bf  2000 AMS Classification Numbers:} 60G15,  60H05, 60F05, 60H07.

 \vskip0.3cm


\section{Introduction}
A famous result by Hsu and Robbins \cite{Hsu} says that if $ X_{1}, X_{2}, \ldots $ is a sequence of independent identically distributed random variables with zero mean and finite variance and $S_{n}:=X_{1}+\ldots +X_{n}$, then
\begin{equation*}
\sum_{n\geq 1} P\left( \vert S_{n}\vert >\eps n\right) <\infty
\end{equation*}
for every $\eps >0$. Later, Erd\"os (\cite{Erd1}, \cite{Erd2}) showed that the converse implication also holds,
namely if the above series is finite for every $\eps>0$ and $X_{1}, X_{2}, \ldots $ are independent and identically
distributed, then $\mathbf{E}X_{1}=0$ and $\mathbf{E}X_{1}^{2}<\infty$. Since then, many authors extended this result
in several directions.

Spitzer's showed in \cite{Spi} that
\begin{equation*}
\sum_{n\geq 1} \frac{1}{n}P\left( \vert S_{n}\vert >\eps n\right) <\infty
\end{equation*}
for every $\eps >0$ if and only if $\mathbf{E}X_{1}=0$ and
$\mathbf{E}\vert X_{1}\vert <\infty. $ Also, Spitzer's theorem has
been the object of various generalizations and variants. One   of the
problems related to the Hsu-Robbins' and Spitzer's theorems is to
find the precise asymptotic as $\eps \to 0$ of the quantities
$\sum_{n\geq 1} P\left( \vert S_{n}\vert >\eps n\right) $ and
$\sum_{n\geq 1} \frac{1}{n}P\left( \vert S_{n}\vert >\eps n\right)$.
Heyde \cite{Hey} showed that
\begin{equation}
\label{heyde}
\lim_{\eps \to 0}\eps^{2}\sum_{n\geq 1} P\left( \vert S_{n}\vert >\eps n\right)=\mathbf{E}X_{1}^{2}
\end{equation}
whenever $\mathbf{E}X_{1}=0$ and $\mathbf{E}X_{1}^{2}<\infty$. In the case when $X$ is attracted to a stable distribution of exponent $\alpha >1$, Spataru \cite{Spa} proved that
\begin{equation}
\label{spataru}
\lim_{\eps \to 0}\frac{1}{-\log \eps } \sum_{n\geq 1} \frac{1}{n}P\left( \vert S_{n}\vert >\eps n\right) =\frac{\alpha}{\alpha-1}.
\end{equation}
The purpose of this paper is to prove Hsu-Robbins and Spitzer's theorems for sequences of correlated random variables, related to the increments of fractional Brownian motion, in the spirit of \cite{Hey} or \cite{Spa}. Recall that the fractional Brownian motion $(B^{H}_{t})_{t\in [0,1]}$ is a centered Gaussian process with covariance function $R^{H}(t,s)=\mathbf{E}(B^{H}_{t}B^{H}_{s} )=\frac{1}{2} (t^{2H}+s^{2H}-\vert t-s\vert ^{2H})$. It can be also defined as the unique self-similar Gaussian process with stationary increments. Concretely, in this paper we will study the behavior of the tail probabilities of the sequence
\begin{equation}\label{vn}
V_{n}= \sum_{k=0}^{n-1} H_{q}\left( n^{H}\left( B_{\frac{k+1}{n}} -B_{\frac{k}{n}} \right) \right)
\end{equation}
 where $B$ is a fractional Brownian motion with Hurst parameter $H\in (0,1)$ (in the sequel we will omit the superscript $H$ for $B$) and $H_{q} $ is the Hermite polynomial of degree $q\geq 1$ given by $H_{q}(x)= (-1) ^{q} e^{\frac{x^{2}}{2}} \frac{d^{q}}{dx^{q}}(e^{-\frac{x^{2}}{2}}).$ The sequence $V_{n}$ behaves as follows (see e.g. \cite{NNT}, Theorem 1; the result is also recalled in Section 3 of our paper): if $0<H<1-\frac{1}{2q}$, a central limit theorem holds for the renormalized sequence $Z^{(1)}_{n}=\frac{V_{n}}{c_{1,q,H}\sqrt{n}}$ while if $1-\frac{1}{2q}<H<1$,  the sequence $Z^{(2)}_{n}=\frac{V_{n}}{c_{2,q,H}n^{1-q(1-H)}}$ converges in $L^{2}(\Omega)$  to a Hermite random variable of order $q$ (see Section 2 for the definition of the Hermite random variable and Section 3 for a rigorous statement concerning the convergence of $V_{n}$). Here $c_{1,q,H}, c_{2,q,H}$ are explicit positive constants depending on $q$ and $H$.

We note that the techniques generally used in the literature to
prove the Hsu-Robbins and Spitzer's results are strongly related to
the independence of the random variables $X_{1}, X_{2}, \ldots.$ In
our case the variables are correlated. Indeed, for any $k,l\geq 1$
we have \\ $\mathbf{E}\left( H_{q}(B_{k+1} -B_{k} )H_{q}(B_{l+1} -B_{l}
)\right) =\frac{1}{(q!)^{2}}\rho_{H}(k-l)$ where the correlation
function is  $\rho_{H}(k)= \frac{1}{2}\left( (k+1)^{2H}
+(k-1)^{2H}-2k^{2H}\right)$ which is not equal to zero unless
$H=\frac{1}{2}$ (which is the case of the standard Brownian motion).
We use new techniques based on the estimates for the multiple
Wiener-It\^o integrals obtained in \cite{BreNo}, \cite{NoPe1} via
Stein's method and Malliavin calculus. Concretely, we study in this paper the behavior as $\eps \to 0$ of the quantities \begin{equation}\label{in1}
 \sum_{n\geq 1} \frac{1}{n} P\left( V_{n} >\eps n\right) =\sum_{n\geq 1} \frac{1}{n} P\left( Z_{n} ^{(1) } >c_{1,q,H}^{-1}\varepsilon \sqrt{n} \right),
 \end{equation}
and
\begin{equation}\label{in2}
\sum_{n\geq 1} P\left( V_{n} >\eps n\right) =\sum_{n\geq 1}  P\left( Z_{n} ^{(1) } >c_{1,q,H}^{-1}\varepsilon \sqrt{n} \right),
\end{equation}
if $0<H<1-\frac{1}{2q}$ and of
\begin{equation}\label{in3}
  \sum_{n\geq 1} \frac{1}{n} P\left( V_{n} >\eps n^{2-2q(1-H)}\right)=\sum_{n\geq 1} \frac{1}{n} P\left( Z_{n} ^{(2) } >c_{2,q,H}^{-1}\varepsilon n^{1-q(1-H)} \right)
\end{equation}
and
\begin{equation}\label{in4}
  \sum_{n\geq 1}  P\left( V_{n} >\eps n^{2-2q(1-H)}\right)=\sum_{n\geq 1} P\left( Z_{n} ^{(2) } >c_{2,q,H}^{-1}\varepsilon n^{1-q(1-H)} \right)
\end{equation}
if $1-\frac{1}{2q}<H<1$. The basic idea in our proofs is that, if we replace $Z^{(1)} _{n}$ and $Z^{(2)}_{n}$ by their limits (standard normal random variable or Hermite random variable) in the above expressions, the behavior as $\eps \to 0$ can be obtained by standard calculations. Then we need to estimate the difference between the tail probabilities of $Z^{(1)}_{n}, Z^{(2)}_{n}$ and the tail probabilities of their limits. To this end, we will use the estimates obtained in \cite{BreNo}, \cite{NoPe1} via  Malliavin calculus and we are able to prove that this difference converges to zero in all cases. We obtain that, as $\eps \to 0$, the quantities (\ref{in1}) and (\ref{in3}) are of order of $\log \eps $ while the functions (\ref{in2}) and (\ref{in4}) are of order of $\eps^{2}$  and $\eps ^{1-q(1-H)}$  respectively.

 The paper is organized as
follows. Section 2 contains some preliminaries on the stochastic analysis  on Wiener
chaos. In Section 3 we prove the Spitzer's theorem for the variations of
the fractional Brownian motion while Section 4 is devoted to  the
Hsu-Robbins theorem for this sequence.

 Throughout the paper we will
denote by $c$ a generic strictly positive constant which may vary
from line to line (and even on the same line).

\section{Preliminaries}

 Let $(W_{t})_{t\in [0,1]}$ be a classical
Wiener process on a standard Wiener space $\left( \Omega
,{\mathcal{F}},\mathbf{P}\right) $. If $f\in
L^{2}([0,1]^{n})$ with $n\geq 1$ integer, we introduce the multiple Wiener-It\^{o} integral of $f$ with respect to $W$. The basic reference is \cite{N}. 

Let $f\in {\mathcal{S}_{m}}$ be an elementary function with $m$
variables  that can be written as \begin{equation*}
f=\sum_{i_{1},\ldots ,i_{m}}c_{i_{1},\ldots i_{m}}1_{A_{i_{1}}\times
\ldots \times A_{i_{m}}}
\end{equation*}%
where the coefficients satisfy $c_{i_{1},\ldots i_{m}}=0$ if two indices $%
i_{k}$ and $i_{l}$ are equal and the sets $A_{i}\in
{\mathcal{B}}([0,1])$ are disjoint. For  such a step function $f$ we
define

\begin{equation*}
I_{m}(f)=\sum_{i_{1},\ldots ,i_{m}}c_{i_{1},\ldots
i_{m}}W(A_{i_{1}})\ldots W(A_{i_{m}})
\end{equation*}%
where we put $W(A)=\int_{0}^{1}1_{A}(s)dW_{s}$ if $A \in {\mathcal{B}}([0,1])$. It can be seen that the mapping $%
I_{n}$ constructed above from ${\mathcal{S}}_{m}$ to $L^{2}(\Omega
)$ is an
isometry on ${\mathcal{S}}_{m}$ , i.e.%
\begin{equation}
\mathbf{E}\left[ I_{n}(f)I_{m}(g)\right] =n!\langle f,g\rangle
_{L^{2}([0,1]^{n})}\mbox{ if }m=n  \label{isom}
\end{equation}%
and%
\begin{equation*}
\mathbf{E}\left[ I_{n}(f)I_{m}(g)\right] =0\mbox{ if }m\not=n.
\end{equation*}%

Since the set ${\mathcal{S}_{n}}$ is dense in $L^{2}([0,1]^{n})$ for every $n\geq 1$ the mapping $%
I_{n}$ can be extended to an isometry from $L^{2}([0,1]^{n})$ to $%
L^{2}(\Omega)$ and the above properties hold true for this
extension.

We will need the following bound for the tail probabilities of multiple Wiener-It\^o integrals (see \cite{Maj}, Theorem 4.1)
\begin{equation}
\label{tail}
P\left( \vert I_{n}(f)\vert >u \right) \leq c\exp \left(  \left( \frac{-cu}{\sigma }\right) ^{\frac{2}{n}}\right)
\end{equation}
for all $u>0$, $n\geq1$, with $\sigma=\Vert f\Vert _{L^{2}([0,1]^{n})}$.

The Hermite random  variable of order $q\geq 1$ that appears as limit in
Theorem \ref{main}, point ii. is defined as (see \cite{NNT})
\begin{equation}
\label{hermite}
Z= d(q,H)I_{q}(L)
\end{equation}
where the kernel $L\in L^{2}([0,1]^{q}) $ is given by
\begin{equation*}
L(y_{1}, \ldots , y_{q})= \int_{y_{1}\vee \ldots \vee y_{q}} ^{1}\partial_{1}K^{H}(u,y_{1})\ldots  \partial_{1}K^{H}(u,y_{q})du.
\end{equation*}
The constant $d(q,H)$ is a positive normalizing constant that guarantees that $\mathbf{E}Z^{2}=1$ and $K^{H}$ is the standard kernel of the fractional Brownian motion (see \cite{N}, Section 5). We will not need the explicit expression of this kernel. Note that the case $q=1$ corresponds to the fractional Brownian motion and the case $q=2$ corresponds to the Rosenblatt process.

\section{Spitzer's theorem}
Let us start by recalling the following result on the convergence of the sequence $V_{n}$ (\ref{vn}) (see \cite{NNT}, Theorem 1).

\begin{theorem}\label{main} Let $q\geq 2$ an integer and let $(B_{t})_{t\geq 0}$ a fractional Brownian motion
 with Hurst parameter $H\in (0,1)$. Then, with some explicit positive constants $c_{1,q,H}, c_{2,q,H}$
 depending only on $q$ and $H$ we have
\begin{description}
\item{i. }If $0<H<1-\frac{1}{2q}$ then
\begin{equation}
\label{c1}
\frac{V_{n}}{c_{1,q,H}\sqrt{n}} {\overset{\mathrm{Law}}{\longrightarrow}} _{n\to \infty} N(0,1)
\end{equation}
\item{ii. } If $1-\frac{1}{2q}<H<1$ then
\begin{equation}
\label{c2}
\frac{V_{n}}{c_{2,q,H}n^{1-q(1-H)}} {\overset{\mathrm{L^{2}}}{\longrightarrow}} _{n\to \infty} Z
\end{equation}
where $Z$ is a Hermite random variable given by (\ref{hermite}).

\end{description}
\end{theorem}
In the case $H=1-\frac{1}{2q}$ the limit is still Gaussian but the normalization is different.
However we will not treat this case in the present work.

We set
\begin{equation}
\label{zi}
Z_{n}^{(1)} = \frac{V_{n}}{c_{1,q,H}\sqrt{n}}, \hskip0.2 cm Z_{n}^{(2)} = \frac{V_{n}}{c_{2,q,H}n^{1-q(1-H)}}
\end{equation}
with the constants $c_{1,q,H}, c_{2,q,H}$ from Theorem \ref{main}.

Let us denote, for every $\varepsilon >0$,
\begin{equation}
\label{f1}
f_{1} (\varepsilon )= \sum_{n\geq 1} \frac{1}{n} P\left( V_{n} >\eps n\right) =\sum_{n\geq 1} \frac{1}{n} P\left( Z_{n} ^{(1) } >c_{1,q,H}^{-1}\varepsilon \sqrt{n} \right)
\end{equation}
and
\begin{equation}
\label{f2}
f_{2}(\eps) =  \sum_{n\geq 1} \frac{1}{n} P\left( V_{n} >\eps n^{2-2q(1-H)}\right)=\sum_{n\geq 1} \frac{1}{n} P\left( Z_{n} ^{(2) } >c_{2,q,H}^{-1}\varepsilon n^{1-q(1-H)} \right)
\end{equation}

\begin{remark}
It is natural to consider the tail probability of order $n^{2-2q(1-H)}$ in (\ref{f2}) because the $L^{2}$ norm of the sequence $V_{n}$ is in this case of order $n^{1-q(1-H)}$.
\end{remark}

 We are interested to study the behavior of $f_{i}(\varepsilon)$ ($i=1,2$) as $\varepsilon \to 0$.  For a given random variable $X$,  we set $\Phi _{X}(z)= 1-P(X<z)+ P(X<-z)$.

 The first lemma gives the asymptotics of the functions $f_{i}(\epsilon)$ as $\eps \to 0$ when $Z_{n}^{(i)}$ are replaced by their limits.

\begin{lemma}\label{spi1}Consider $c>0$.
\begin{description}
\item{i. }Let $Z^{(1)}$ be  a standard normal random variable. Then as
\begin{equation*}
\frac{1}{-\log c\varepsilon} \sum_{n\geq 1} \frac{1}{n}\Phi _{Z^{(1)}} (c\varepsilon \sqrt{n} ) \to _{\varepsilon \to 0}2.
\end{equation*}
\item{ii. } Let $Z^{(2)}$ be a Hermite random variable or order $q$ given by (\ref{hermite}). Then, for any integer $q\geq 1$
\begin{equation*}
\frac{1}{-\log c\varepsilon} \sum_{n\geq 1} \frac{1}{n}\Phi _{Z^{(2)}} (c\varepsilon n^{1-q(1-H)} ) \to _{\varepsilon \to 0}\frac{1}{1-q(1-H)}.
\end{equation*}
\end{description}
\end{lemma}
{\bf Proof: } The case when $Z^{(1)}$ follows the standard normal
law is hidden in  \cite{Spa}. We will give the ideas of the proof.
We can write (see \cite{Spa})
\begin{equation*}
\sum_{n\geq1}  \frac{1}{n} \Phi _{Z^{(1)}} (c\varepsilon  \sqrt{n} )=\int_{1}^{\infty} \frac{1}{x} \Phi_{Z^{(1)}} (c\varepsilon \sqrt{x} )dx-\frac{1}{2} \Phi _{Z^{(1)}} (c\varepsilon) -\int_{1}^{\infty} P_{1}(x) d\left[ \frac{1}{x} \Phi_{Z^{(1)}} (c\varepsilon \sqrt{x})\right].
\end{equation*}
with $P_{1}(x)=[x]-x+\frac{1}{2}$. Clearly as $\varepsilon \to 0$, $\frac{1}{\log \varepsilon } \Phi _{Z^{(1)}} (c\varepsilon)\to 0$ because $\Phi _{Z^{(1)}}$ is a bounded function and concerning the last term it is also trivial to observe that
\begin{eqnarray*}
&&\frac{1}{-\log c\varepsilon } \int_{1}^{\infty} P_{1}(x) d\left[ \frac{1}{x} \Phi_{Z^{(1)}} (c\varepsilon \sqrt{x})\right]\\
&&=\frac{1}{-\log c\varepsilon } \left(-\int_{1}^{\infty} P_{1}(x) \left(\frac{1}{x^{2}} \Phi_{Z^{(1)}}
 (c\varepsilon \sqrt{x}) dx + c\varepsilon \frac{1}{2}x^{-\frac{1}{2}}
 \frac{1}{x}\Phi_{Z^{(1)}}'(\varepsilon \sqrt{x})\right) dx\right) \to _{\eps \to 0} 0
\end{eqnarray*}
since $\Phi _{Z^{(1)}} $ and $\Phi_{Z^{(1)}}' $ are bounded.  Therefore the asymptotics of the function $f_{1}(\eps)$ as  $\eps \to 0$ will be given by $\int_{1}^{\infty} \frac{1}{x} \Phi_{Z^{(1)}} (c\eps \sqrt{x}) dx.$
By making the change of variables $c\eps \sqrt{x}=y$, we get
\begin{eqnarray*}
&&\lim _{\eps \to 0} \frac{1}{-\log c\eps }\int_{1}^{\infty} \frac{1}{x} \Phi_{Z^{(1)}} (c\eps \sqrt{x}) dx
= \lim _{\eps \to 0} \frac{1}{-\log c\eps }2\int_{c\eps }^{\infty } \frac{1}{y} \Phi_{Z^{(1)}} (y) dy =\lim _{\eps \to 0} 2\Phi_{Z^{(1)}}(c\eps)=2.
\end{eqnarray*}

\vskip0.2cm

Let us consider now  the case of the Hermite random variable. We will have as above
 \begin{eqnarray*}
 &&\lim_{\eps \to 0} \frac{1}{-\log c\eps }\sum_{n\geq 1}\frac{1}{n}\Phi_{Z^{(2)}}(c\eps n^{1-q(1-H)} ) \\
 &=& \lim_{\eps \to 0} \frac{1}{-\log c\eps } \left( \int_{1}^{\infty} \frac{1}{x} \Phi_{Z^{(2)}}( c\eps x^{1-q(1-H)} )dx -\int_{1}^{\infty} P_{1}(x) d\left[ \frac{1}{x} \Phi_{Z^{(2)}}( c\eps x^{1-q(1-H)} )\right]\right)
 \end{eqnarray*}
By making the change of variables $c\eps x^{1-q(1-H)} =y$ we will obtain
\begin{eqnarray*}
&&\lim_{\eps \to 0} \frac{1}{-\log c\eps }  \int_{1}^{\infty} \frac{1}{x} \Phi_{Z^{(2)}}( c\eps x^{1-q(1-H)} )dx\\
&=& \lim_{\eps \to 0} \frac{1}{-\log c\eps }\frac{1}{1-q(1-H)}\int_{c\eps } ^{\infty } \frac{1}{y} \Phi_{Z^{(2)}} (y) dy = \lim_{\eps \to 0}\frac{1}{1-q(1-H)} \Phi_{Z^{(2)}}(c\eps ) = \frac{1}{1-q(1-H)}
\end{eqnarray*}
where we used the fact that $\Phi_{Z^{(2)}}(y) \leq y^{-2} \mathbf{E}\vert Z^{(2)}\vert ^{2} $ and so $ \lim_{y\to \infty} \log y \Phi_{Z^{(2)}}(y) =0$.

It remains to show that
$\frac{1}{-\log c\eps } \int_{1}^{\infty} P_{1}(x) d\left[ \frac{1}{x} \Phi_{Z^{(2)}}( c\eps x^{1-q(1-H)} )\right]$
converges to zero as $\eps $ tends to 0 (note that actually it follows from a result by \cite{Albin} that a Hermite
 random variable has a density, but we don't need it explicitly, we only use the fact that $\Phi_{Z^{(2)}}$
 is differentiable almost everywhere). This  is equal to
\begin{eqnarray*}
 &&\lim_{\eps }\frac{1}{-\log c\eps }\int_{1}^{\infty} P_{1}(x) c\eps (1-q(1-H)) x^{-q(1-H)-1} \Phi_{Z^{(2)}}'( c\eps x^{1-q(1-H)} )dx\\
&=& c\frac{\eps }{-\log \eps } (c\eps) ^{\frac{q(1-H)}{1-q(1-H)}} \int_{c\eps }^{\infty} P_{1}\left(  \left( \frac{y}{c\eps}\right) ^{\frac{1}{1-q(1-H)}}\right) \Phi_{Z^{(2)}}'(y) y^{-\frac{1}{1-q(1-H)}} dy \\
&\leq & c\frac{1}{-\log \eps }\int_{c \eps}^{\infty} P_{1}\left(  \left( \frac{1}{c\eps}\right) ^{\frac{1}{1-q(1-H)}}\right) \Phi_{Z^{(2)}}'(y) dy
\end{eqnarray*}
which clearly goes to zero since $P_{1}$ is bounded and $\int_{0}^{\infty} \Phi _{Z^{(2)}}'(y) dy=1$.\qed

\vskip0.3cm

The next result estimates the limit of the difference between the functions $f_{i}(\eps)$ given by (\ref{f1}), (\ref{f2}) and the sequence in Lemma 1.
 \begin{prop}\label{spi2}Let $q\geq 2$ and $c>0$.
 \begin{description}
 \item{i. } If $H<1-\frac{1}{2q}$, let $Z_{n}^{(1)}$ be given by (\ref{zi}) and let $Z^{(1)}$ be standard normal random variable. Then it holds
 \begin{equation*}
 \frac{1}{-\log c\varepsilon } \left[    \sum_{n\geq 1}\frac{1}{n} P\left( \vert Z_ {n}^{(1)}\vert  > c\varepsilon \sqrt{n}\right) - \sum_{n\geq 1}\frac{1}{n} P\left( \vert Z ^{(1)}\vert  > c\varepsilon \sqrt{n} \right)\right]\to_{\eps\to 0}0.
 \end{equation*}
 \item{ii. } Let $Z^{(2)}$ be a Hermite random variable of order $q \geq 2$ and $H>1-\frac{1}{2q}$. Then
 \begin{equation*}
 \frac{1}{-\log c\varepsilon } \left[    \sum_{n\geq 1}\frac{1}{n} P\left( \vert Z_ {n}^{(2)}\vert  > c\varepsilon n^{1-q(1-H)} \right) - \sum_{n\geq 1}\frac{1}{n} P\left( \vert Z^{(2)}\vert  > c\varepsilon n^{1-q(1-H)} \right)\right]\to_{\eps\to 0}0.
 \end{equation*}
 \end{description}

 \end{prop}
 {\bf Proof: } Let us start with the point i. Assume  $H<1- \frac{1}{2q}.$ We can write
 \begin{eqnarray*}
 && \sum_{n\geq 1}\frac{1}{n} P\left( \vert Z_{n} ^{(1)} \vert > c\varepsilon \sqrt{n} \right) - \sum_{n\geq 1}\frac{1}{n} P\left( \vert Z ^{(1)}\vert  > c\varepsilon \sqrt{n} \right)\\
 &=&\sum_{n\geq 1} \frac{1}{n}\left[ P\left(  Z_{n} ^{(1)} > c\varepsilon \sqrt{n} \right) -  P\left(  Z ^{(1)}  > c\varepsilon \sqrt{n} \right)\right]
 + \sum_{n\geq 1}\left[\frac{1}{n} P\left(  Z_{n} ^{(1)}  <- c\varepsilon \sqrt{n} \right) -  P\left(  Z ^{(1)}<- c\varepsilon \sqrt{n}\right)\right]\\
 &\leq & 2\sum_{n\geq 1} \frac{1}{n} \sup_{x\in \mathbb{R}}\left|   P\left(  Z_{n} ^{(1)} >x \right)- P\left(  Z ^{(1)} >x \right)\right|.
 \end{eqnarray*}

 It follows from \cite{NoPe1}, Theorem 4.1 that
 \begin{equation}\label{bound1}
 \sup_{x\in \mathbb{R}}\left|   P\left(  Z_{n} ^{(1)} >x \right)- P\left(  Z ^{(1)} >x \right)\right|\leq c\left\{
  \begin{array}{lcl}
\frac{1}{\sqrt{n}}, \hskip0.3cm  H\in(0, \frac{1}{2}] \\
n^{H-1}, \hskip0.3cm H\in [\frac{1}{2}, \frac{2q-3}{2q-2}) \\
n^{qH-q+\frac{1}{2} }, \hskip0.3cm H\in [\frac{2q-3}{2q-2}, 1-\frac{1}{2q}).
  \end{array}\right.
 \end{equation}
  and this implies that
  \begin{equation}
  \sum_{n\geq 1} \frac{1}{n} \sup_{x\in \mathbb{R}}\left|   P\left(  Z_{n} ^{(i)} >x \right)- P\left(  Z ^{(i)} >x \right)\right|
 \leq c \left\{
  \begin{array}{lcl}
\sum_{n\geq 1} \frac{1}{n\sqrt{n}}, \hskip0.3cm  H\in(0, \frac{1}{2}] \\
\sum_{n\geq 1}n^{H-2}, \hskip0.3cm H\in [\frac{1}{2}, \frac{2q-3}{2q-2}) \\
\sum_{n\geq 1}n^{qH-q-\frac{1}{2} }, \hskip0.3cm H\in [\frac{2q-3}{2q-2}, 1-\frac{1}{2q}).
  \end{array}\right.
 \end{equation}
 and the last sums are finite (for the last one we use $H<1-\frac{1}{2q}$). The conclusion follows.

 \vskip0.3cm

 Concerning the point ii. (the case $H>1-\frac{1}{2q}$), by using  a result in Proposition 3.1 of \cite{BreNo} we have
 \begin{equation}\label{bound2}
 \sup_{x\in \mathbb{R}}\left|   P\left(  Z_{n} ^{(i)} >x \right)- P\left(  Z ^{(i)} >x \right)\right|\leq c\left( \mathbf{E}\left| Z_{n}^{(2)}-Z^{(2)}\right| ^{2}\right) ^{\frac{1}{2q}}\leq cn^{1-\frac{1}{2q}-H}
 \end{equation}
 and as a consequence
 \begin{equation*}
  \sum_{n\geq 1}\frac{1}{n} P\left( \vert Z_{n} ^{(2)} \vert > c\varepsilon n^{1-q(1-H)} \right) - \sum_{n\geq 1}\frac{1}{n} P\left( \vert Z ^{(2)}\vert  > c\varepsilon n^{1-q(1-H)} \right) \leq c \sum_{n\geq 1} n^{-\frac{1}{2q}-H}
  \end{equation*}
and the above series is convergent because $H>1-\frac{1}{2q}$. \qed

\vskip0.5cm
We state now the Spitzer's theorem for the variations of the fractional Brownian motion.
\begin{theorem} Let  $f_{1}, f_{2}$ be given by (\ref{f1}), (\ref{f2}) and the constants $c_{1,q,H}, c_{2,q,H}$ be those from Theorem \ref{main}.
\begin{description}
\item{i. }If $0<H<1-\frac{1}{2q}$ then
\begin{equation*}
\lim_{\eps \to 0} \frac{1}{\log (c_{1,H,q}^{-1}\eps)} f_{1}(\eps )=2.
\end{equation*}
\item{ii. }If  $1>H>1-\frac{1}{2q}$  then
\begin{equation*}
\lim_{\eps \to 0} \frac{1}{\log (c_{2,H,q}^{-1}\eps)} f_{2}(\eps )=\frac{1}{1-q(1-H)}.
\end{equation*}
\end{description}

\end{theorem}
{\bf Proof: } It is a consequence of Lemma \ref{spi1} and Proposition \ref{spi2}. \qed

\vskip0.2cm

\begin{remark}
Concerning the case $H=1-\frac{1}{2q}$, note that the correct
normalization of $V_{n}$ (\ref{vn}) is $\frac{1}{(\log n)\sqrt{n}}$. Because of
the appearance  of the term $\log n$  our approach is not directly applicable to this case.
\end{remark}

 \section{Hsu-Robbins theorem for the variations of  fractional Brownian motion}

In this section we prove a version of the Hsu-Robbins theorem for the variations of the fractional Brownian motion.  Concretely, we denote here by, for every $\eps >0$
 \begin{equation}\label{g1}
 g_{1}(\varepsilon )= \sum_{n\geq 1} P \left( \vert V_{n}\vert  > \eps n\right)
 \end{equation}
 if $H<1-\frac{1}{2q}$ and by
 \begin{equation}\label{g2}
 g_{2}(\varepsilon )= \sum_{n\geq 1} P \left( \vert V_{n}\vert  > \eps n^{2-2q(1-H)}\right)
 \end{equation}
if $H>1-\frac{1}{2q}$.
and we  estimate the behavior of the functions $g_{i}(\varepsilon) $ as $\eps \to 0$. Note that we can write
 \begin{equation*}
 g_{1}(\eps) = \sum_{n\geq 1} P\left( \vert Z_{n}^{(1)}\vert  > c_{1,q,H}^{-1} \eps \sqrt{n}\right), \hskip0.2cm g_{2}(\eps) =  \sum_{n\geq 1} P\left( \vert Z_{n}^{(2)}\vert  > c_{2,q,H}^{-1} \eps n^{1-q(1-H)}\right)
 \end{equation*}
 with $Z_{n}^{(1)}, Z_{n}^{(2)}$ given by (\ref{zi}).

 We decompose it as: for $H<1-\frac{1}{2q}$
 \begin{eqnarray*}
 g_{1}(\eps )&=& \sum_{n\geq 1}  P \left( \vert Z^{(1)} \vert >c_{1,q,H}^{-1} \eps \sqrt {n}\right)\\
 &+& \sum_{n\geq 1} \left[ P \left( \vert Z_{n}^{(1)} \vert > c_{1,q,H}^{-1}\eps \sqrt{n}\right) -P \left( \vert Z^{(1)}\vert  >c_{1,q,H}^{-1} \eps \sqrt{n}\right)\right].
 \end{eqnarray*}
 and for $H>1-\frac{1}{2q}$
 \begin{eqnarray*}
 g_{2}(\eps )&=& \sum_{n\geq 1}  P \left( \vert Z^{(2)} \vert > \eps c_{2,q,H}^{-1}n^{1-q(1-H)}\right)\\
 &+& \sum_{n\geq 1} \left[ P \left( \vert Z_{n}^{(2)} \vert >c _{2,q, H}^{-1}\eps n^{1-q(1-H)}\right) -P \left( \vert Z^{(2)}\vert  > c_{2,q, H}^{-1}\eps n^{1-q(1-H)}\right)\right].
 \end{eqnarray*}
We start again by consider the situation when $Z^{(i)}_{n}$ are
replaced by their limits.
 \begin{lemma}\label{hr1}
 \begin{description}
 \item{i. }
 Let $Z^{(1)} $ be a standard normal random variable. Then
 \begin{equation*}
\lim_{\eps \to 0} (c\eps) ^{2} \sum_{n\geq 1}  P \left( \vert Z^{(1)} \vert > c\eps \sqrt{n} \right) =1.
 \end{equation*}
 \item{ii. } Let $Z^{(2)}$ be a Hermite random variable with $H>1-\frac{1}{2q}$.  Then
 \begin{equation*}
 \lim_{\eps \to 0} (c\eps ) ^{\frac{1}{1-q(1-H)}} \sum _{n\geq 1} P\left( \vert Z^{(2)}\vert > c \eps n^{1-q(1-H)} \right) = \mathbf{E}\vert Z^{(2)}\vert ^{\frac{1}{1-q(1-H)}}.
 \end{equation*}
 \end{description}
 \end{lemma}
 {\bf Proof: } The part i.  is a consequence of the result of Heyde \cite{Hey}. Indeed take $X_{i}\sim N(0,1)$ in (\ref{heyde}).  Concerning part ii. we can write
 \begin{eqnarray*}
  &&\lim_{\eps \to 0} (c\eps ) ^{\frac{1}{1-q(1-H)}} \sum_{n\geq 1} \Phi_{Z^{(2)}} (c\eps n^{1-q(1-H)}) \\
  &=& \lim_{\eps \to 0} (c\eps ) ^{\frac{1}{1-q(1-H)}}\left[ \int_{1}^{\infty} \Phi_{Z^{(2)}} (c\eps x^{1-q(1-H)})dx -\int_{1}^{\infty} P_{1}(x) d\left[ \Phi_{Z^{(2)}}(c\eps x^{1-q(1-H)})\right]\right] \\
  &&:=\lim_{\eps \to 0} (A(\eps) + B(\eps))
  \end{eqnarray*}
 with $P_{1}(x)= [x]-x+\frac{1}{2}$. Moreover
 \begin{eqnarray*}
 A(\eps)&=& (c\eps ) ^{\frac{1}{1-q(1-H)}} \int_{1}^{\infty} \Phi_{Z^{(2)}} (c\eps x^{1-q(1-H)})dx\\
 &=&\frac{1}{1-q(1-H)} \int_{c\eps }^{\infty} \Phi_{Z^{(2)}}(y) y^{\frac{1}{1-q(1-H)} -1}dy .
\end{eqnarray*}
 Since $\Phi_{Z^{(2)}}(y)\leq y^{-2}$ we have $\Phi_{Z^{(2)}}(y)y^{\frac{1}{1-q(1-H)} }\to_{y\to \infty}0$ and therefore
 $$A(\eps) = -\Phi_{Z^{(2)}}(c\eps) (c\eps) ^{\frac{1}{1-q(1-H)}}- \int_{c\eps }^{\infty} \Phi_{Z^{(2)}}'(y) y^{\frac{1}{1-q(1-H)} }dy$$
 where the first terms goes to zero and the second to $\mathbf{E}\left| Z^{(2)}\right|  ^{\frac{1}{1-q(1-H)}}.$ The proof that the term $B(\eps)$ converges to zero is similar to the proof of Lemma \ref{hr1}, point ii.\qed

\vskip0.3cm

\begin{remark}
 The Hermite random variable has moments of all orders (in particular the moment of order $\frac{1}{1-q(1-H)}$ exists) since it is the value at time 1  of a selfsimilar process with stationary increments.
\end{remark}

 \begin{prop}\label{hr2}
 \begin{description}
 \item{i. }Let $H<1-\frac{1}{2q}$ and let $Z_{n}^{(1)}$ be given by (\ref{zi}). Let also $Z^{(1)}$ be a standard normal random variable. Then
 \begin{equation*}(c\eps )^{2}\sum_{n\geq 1} \left[ P \left( \vert Z_{n}^{(1)} \vert > c\eps \sqrt{n}\right) -P \left( \vert Z^{(1)}\vert  >c \eps \sqrt{n}\right)\right]\to _{\eps \to 0}0
  \end{equation*}
 \item{ii. } Let $H>1-\frac{1}{2q}$ and let $Z_{n}^{(2)}$ be given by (\ref{zi}). Let $Z^{(2)}$ be a Hermite random variable. Then
 \begin{equation*}
 (c\eps ) ^{\frac{1}{1-q(1-H)}}
 \sum_{n\geq 1} \left[ P \left( \vert Z_{n}^{(2)} \vert > c\eps n^{1-q(1-H)}\right) -P \left( \vert Z^{(2)}\vert  >c \eps n^{1-q(1-H)}\right)\right]\to _{\eps \to 0} 0.
 \end{equation*}
 \end{description}
 \end{prop}

 \begin{remark}
 Note that the bounds (\ref{bound1}), (\ref{bound2}) does not help here because the series that appear after their use are not convergent.
\end{remark}

\b {\bf Proof of Proposition \ref{hr2}: } {\it Case
$H<1-\frac{1}{2q}$. }We have, for some $\beta >0$ to be chosen
later,
\begin{eqnarray*}
&&\eps ^{2}
 \sum_{n\geq 1} \left[ P \left( \vert Z_{n}^{(1)} \vert > c\eps \sqrt{n}\right) -P \left( \vert Z^{(1)}\vert  >c \eps \sqrt{n}\right)\right]\\
 &=&
 \eps ^{2}
 \sum_{n=1}^{[\eps ^{-\beta }]} \left[ P \left( \vert Z_{n}^{(1)} \vert > c\eps \sqrt{n}\right) -P \left( \vert Z^{(1)}\vert  >c \eps \sqrt{n}\right)\right]\\
 &&+ \eps ^{2}
 \sum_{n> [\eps ^{-\beta}]} \left[ P \left( \vert Z_{n}^{(1)} \vert > c\eps\sqrt{n}\right) -P \left( \vert Z^{(1)}\vert  >c \eps \sqrt{n}\right)\right]\\
 &:=& I_{1}(\eps)+ J_{1}(\eps).
\end{eqnarray*}

  Consider first the situation when $H\in (0, \frac{1}{2}]$. Let us choose a real number $\beta $ such that $2<\beta <4 $.
By using (\ref{bound1}),
\begin{equation*}I_{1}(\eps)\leq c\eps ^{2}  \sum_{n=1}^{[\eps ^{-\beta }]}n^{-\frac{1}{2}} \leq c\eps^{2} \eps^{-\frac{\beta}{2}}\to _{\eps \to 0}0
\end{equation*}
since $\beta <4$.  Next, by using the bound for the tail probabilities of multiple integrals and since $\mathbf{E}\left| Z_{n}^{(1)}\right| ^{2}$ converges to $1$ as $n\to \infty$
\begin{eqnarray*}
&&J_{1}(\eps)=\eps ^{2}  \sum_{n> [\eps ^{-\beta}]}P\left( Z_{n}^{(1)} >c\eps \sqrt{n}\right) \leq c \eps ^{-2}  \sum_{n> [\eps ^{-\beta}]}\exp \left(\frac {-c\eps \sqrt{n}}{\left(\mathbf{E}\left|Z_{n}^{(1)}\right| ^{2} \right)^{\frac{1}{2}}}\right)  ^{\frac{2}{q}}\\
&&\leq \eps ^{2}  \sum_{n> [\eps ^{-\beta}]} \exp \left( \left( -cn^{-\frac{1}{\beta} }\sqrt{n}\right) ^{\frac{2}{q}}\right)
\end{eqnarray*}
and since converges to zero for $\beta >2$. The same argument shows that $\eps ^{2}  \sum_{n> [\eps ^{-\beta}]}P\left( Z^{(1)} >c\eps \sqrt{n}\right)$ converges to zero.

The case when $H\in (\frac{1}{2}, \frac{2q-3}{2q-2}) $ can be
obtained by taking $2<\beta < \frac{2}{H}$ (it is possible since
$H<1$) while in the case $H\in (\frac{2q-3}{2q-2}, 1-\frac{1}{2q})$
we have to choose $2<\beta < \frac{2}{qH-q+\frac{3}{2}}$
 (which is possible because $H<1-\frac{1}{2q}$!).

\vskip0.2cm

\b {\it Case $H>1-\frac{1}{2q}$. }  We have, with some suitable $\beta>0$
\begin{eqnarray*}
&&\eps ^{\frac{1}{1-q(1-H)}}
 \sum_{n\geq 1} \left[ P \left( \vert Z_{n}^{(2)} \vert > c\eps n^{1-q(1-H)}\right) -P \left( \vert Z^{(2)}\vert  >c \eps n^{1-q(1-H)}\right)\right]\\
 &=&
 \eps ^{\frac{1}{1-q(1-H)}}
 \sum_{n=1}^{[\eps ^{-\beta }]} \left[ P \left( \vert Z_{n}^{(2)} \vert > c\eps n^{1-q(1-H)}\right) -P \left( \vert Z^{(2)}\vert  >c \eps n^{1-q(1-H)}\right)\right]\\
 &&+ \eps ^{\frac{1}{1-q(1-H)}}
 \sum_{n\geq [\eps ^{-\beta}]} \left[ P \left( \vert Z_{n}^{(2)} \vert > c\eps n^{1-q(1-H)}\right) -P \left( \vert Z^{(2)}\vert  >c \eps n^{1-q(1-H)}\right)\right]\\
 &:=& I_{2}(\eps)+ J_{2}(\eps).
\end{eqnarray*}
Choose $\frac{1}{1-q(1-H)}<\beta <
\frac{1}{(1-q(1-H))(2-H-\frac{1}{2q})}$ (again, this is always
possible when $H>1-\frac{1}{2q}$!). Then
\begin{equation*}
I_{2}(\eps)\leq c e^{\frac{1}{1-q(1-H)}} \eps^{(-\beta)(2-H-\frac{1}{2q})}\to _{\eps \to 0}0
\end{equation*}
and by (\ref{tail})
\begin{equation*}
J_{2}(\eps)\leq c \sum_{n>[\eps ^{-\beta}]} \exp \left( \left(
\frac{-c\eps n^{1-q(1-H)}}{\left( \mathbf{E}\left|Z_{n}^{(2)}\right|
^{2} \right)^{\frac{1}{2}}}\right) ^{\frac{2}{q}}\right) \leq c
\sum_{n>[\eps ^{-\beta}]} \exp \left(
cn^{-\frac{1}{\beta}}n^{1-q(1-H)}\right) ^{\frac{2}{q}}\to _{\eps
\to 0} 0
\end{equation*}\qed

\vskip0.3cm We state the main result of this section which is a
consequence of  Lemma \ref{hr1} and Proposition \ref{hr2}.
\begin{theorem}
Let $q\geq 2$ and let $c_{1,q,H}, c_{2,q,H}$ be the constants from
Theorem \ref{main}. Let $Z^{(1)} $ be a standard normal random
variable, $Z^{(2)} $ a Hermite random variable of order $q\geq 2$
and let $g_{1}, g_{2}$ be given by (\ref{g1}) and (\ref{g2}). Then
\begin{description}
\item{i. } If $0<H<1-\frac{1}{2q}$, we have $(c_ {1,q,H}^{-1}\eps) ^{2} g_{1}(\eps) \to _{\eps \to 0}1=\mathbf{E}Z^{(1)}$.
\item{ii. } If $1-\frac{1}{2q}<H<1$ we have $(c_{2,q,H}^{-1}\eps) ^{\frac{1}{1-q(1-H)}} g_{2}(\eps)\to _{\eps \to 0} \mathbf{E}\vert Z^{(2)}\vert ^{\frac{1}{1-q(1-H)}}$.
\end{description}
\end{theorem}

\begin{remark}
In the case $H=\frac{1}{2}$ we retrieve the result (\ref{heyde}) of \cite{Hey}. The case $q=1$ is trivial, because  in this case, since $V_{n}=B_{n}$ and $\mathbf{E}V_{n}^{2}=n^{2H}$, we obtain the following (by applying Lemma \ref{spi1} and \ref{hr1} with $q=1$)
\begin{equation*}
\frac{1}{\log \eps} \sum_{n\geq 1} \frac{1}{n}P\left( \vert V_{n}\vert >\eps n^{2H} \right)  \to_{\eps \to 0}\frac{1}{H}
\end{equation*}
and
\begin{equation*}
\eps ^{2}\sum_{n\geq 1} P\left( \vert V_{n}\vert >\eps n^{2H} \right)  \to_{\eps \to 0}\mathbf{E}\left| Z^{(1)} \right| ^{\frac{1}{H}}.
\end{equation*}
\end{remark}

\begin{remark}
Let $(\eps _{i})_{i\in \mathbb{Z}}$ be a sequence of i.i.d. centered random variable with finite variance and let $(a_{i})_{i\geq 1}$ a square summable real sequence. Define $X_{n}=\sum_{i\geq 1}a_{i}\eps _{n-i}$. Then the sequence $S_{N}=\sum_{n=1}^{N} \left[ K(X_{n})-\mathbf{E}K(X_{n})\right]$ satisfies a central limit theorem or a non-central limit theorem according to the properties of the measurable function $K$ (see \cite{HH} or \cite{Wu}). We think that our tools can be applied to investigate the tail probabilities  of the sequence $S_{N}$  in the spirit of \cite{Hey} or \cite{Spa} at least the in particular cases (for example, when $\eps_{i}$ represents the increment $W_{i+1}-W_{i}$ of a Wiener process because in this case $\eps _{i}$ can be written as a multiple integral of order one and $X_{n}$ can be decomposed into a sum of multiple integrals. We thank the referee for mentioning   the references \cite{HH} and \cite{Wu}.
\end{remark}

\end{document}